\documentclass[reqno,centertags,11pt]{amsart}
\usepackage{verbatim}
\usepackage{tikz, enumerate, array, amssymb, amsmath}
\usepackage{graphicx, epsfig }
\usepackage{amsthm}
\usepackage{amsbsy}
\usepackage{amsfonts}

\usepackage[colorlinks]{hyperref}

\usepackage{mathtools}

\newcommand{\Hmm}[1]{\leavevmode{\marginpar{\tiny%
$\hbox to 0mm{\hspace*{-0.5mm}$\leftarrow$\hss}%
\vcenter{\vrule depth 0.1mm height 0.1mm width \the\marginparwidth}%
\hbox to 0mm{\hss$\rightarrow$\hspace*{-0.5mm}}$\\\relax\raggedright #1}}}

\newcommand{\N}{{\mathbb{N}}}

\newcommand{\R}{{\mathbb{R}}}

\newcommand{\C}{{\mathbb{C}}}

\newcommand{\red}{\color{red} }

\newcommand{\f}{\frac}

\newcommand{\beq}{\begin{equation}}
\newcommand{\eeq}{\end{equation}}
\newcommand{\bdm}{\begin{displaymath}}
\newcommand{\edm}{\end{displaymath}}
\newcommand{\ba}{\begin{align}}
\newcommand{\ea}{\end{align}}
\newcommand{\bpf}{\begin{proof}}
\newcommand{\epf}{\end{proof}}

\newcommand{\veps}{\varepsilon}

\newcommand{\im}{\mathrm{Im}}

\newcommand{\dav}{{d_{\mathrm{av}}}}

\allowdisplaybreaks

\newcommand{\calC}{\mathcal{C}}

\newtheorem{theorem}{Theorem}
\newtheorem{proposition}[theorem]{Proposition}
\newtheorem{lemma}[theorem]{Lemma}

\theoremstyle{definition}

\newtheorem{remark}[theorem]{Remark}
\newtheorem{remarks}[theorem]{Remarks}

\newcounter{theoremi}[theorem]

\numberwithin{theorem}{section}
\numberwithin{equation}{section}

\newcounter{smalllist}

\newcounter{listi}
\newenvironment{theoremlist}{\begin{list}{{\rm(\roman{listi})}}{%
\setlength{\topsep}{0mm}\setlength{\parsep}{0mm}\setlength{\itemsep}{0mm}%
\setlength{\labelwidth}{1.5em}\setlength{\leftmargin}{1.7em}\usecounter{listi}%
}}{\end{list}}

\newcounter{smallenum}

\newcounter{assumptions}

\title[Global existence versus finite time blowup dichotomy for DMNLS]{Global existence versus finite time blowup dichotomy for the dispersion managed NLS}

\date{\today}
\author[M. Choi]{Mi-Ran Choi}
\address{Department of Mathematics, Sogang University, 35 Baekbeom-ro,
  Mapo-gu, Seoul 04107, South Korea.}%
\email{mrchoi@sogang.ac.kr}

\author[Y. Hong]{Younghun Hong}
\address{Department of Mathematics, Chung-Ang University, 84 Heukseok-ro, Dongjak-gu, Seoul 06974, South Korea.}
\email{yhhong@cau.ac.kr}

\author[Y-R, Lee]{Young-Ran Lee}
\address{Department of Mathematics, Sogang University, 35 Baekbeom-ro,
  Mapo-gu, Seoul 04107, South Korea.}%
\email{younglee@sogang.ac.kr}

\begin{document}
\maketitle

\begin{abstract}
We consider
  the Gabitov-Turitsyn equation or the dispersion managed nonlinear
  Schr\"odinger equation of a power-type nonlinearity
  \[
   i\partial_t u+\dav \partial_x^2u+\int_0^1 e^{-ir\partial_x^2}\big(|e^{ir\partial_x^2}u|^{p-1}e^{ir\partial_x^2}u\big)dr=0
  \]
and prove the global existence versus finite time blowup dichotomy for the mass-supercritical cases, that is, $p>9$.
\end{abstract}
\section{Introduction}
 In this paper, we consider the Cauchy problem for the dispersion managed nonlinear Schr\"odinger equation (NLS)
 \begin{equation}\label{NLS}
\begin{cases}
        \ \displaystyle{ i\partial_t u+\dav \partial_x^2u+\int_0^1 e^{-ir\partial_x^2}\big(|e^{ir\partial_x^2}u|^{p-1}e^{ir\partial_x^2}u\big)dr=0, }\\[.5ex]
            \  u(0,\cdot) = u_0
            \end{cases}
\end{equation}
where $u=u(t,x): I\times\mathbb{R}\to\mathbb{C}$ for an interval $I$, $p>1$, $\dav\in \R$, and $e^{ir\partial_x^2}$ is the linear propagator.

The model \eqref{NLS} arises from the study of NLS with a periodically varying dispersion coefficient
\beq\label{eq:NLS}
i\partial_t u+ d(t) \partial_x^2u+ |u|^{p-1}u=0.
\eeq
This equation describes the propagation of signals through glass-fiber cables with alternating sections of strongly positive and strongly negative dispersion, the so-called strong dispersion management. Here, $t$ corresponds to the distance along the fiber and $x$ denotes the (retarded) time. Hence, $d(t)$ is not varying in time but represents a dispersion varying along the optical cable.
Specifically, the local dispersion $d(t)$ is given by
 \[
 d(t)=\dav+\veps^{-1}d_0(t/\veps)
 \]
 where $\dav$ is the average component, $d_0$ is its mean zero part over one period and $\veps$ is a small parameter.
The technique of dispersion management was introduced in  \cite{1980} and proved to be incredibly successful in producing stable, soliton-like pulses.
See the reviews \cite{review1, review2} and the references cited in \cite{HL} for a discussion of the dispersion management technique.
 It is by now well-known that equation \eqref{eq:NLS} with strong dispersion management can be averaged over one period to yield an effective equation, see \cite{GT1, GT2}.
In particular, if $d_0$ is the $2$-periodic function with $d_0(t)=1$ on $[0,1)$ and $d_0(t)=-1$ on $[1,2)$, then the change of variables and the averaging process yield the equation of the form \eqref{NLS},
called the Gabitov-Turitsyn equation. This averaging process is verified in \cite{CLA, ZGJT01}.

For the Cauchy problem \eqref{NLS}, it is shown to be locally well-posed in $H^1(\mathbb{R})$ for all $p>1$ when  $\dav\neq0$ ; $1<p \le 5$ when $\dav=0$, see  \cite{CHL}.
Furthermore, its solution preserves the mass and the energy, that is, for all $t$,
\[
M[u(t)]=M[u_0 ] \quad\mbox{and}\quad E[u(t)]=E[u_0 ],
\]
where the mass and energy are defined by
\[
M[f]:=\|f\|_{L^2}^2 \quad\mbox{and}\quad
E[f]:=\f{\dav}{2}\|\partial_x f\|^2_{L^2} - \f{1}{p+1} \int_0^1 \|e^{ir\partial_x^2} f\|^{p+1}_{L^{p+1}}dr.
\]
When $\dav$ is negative, by the conservation laws, solutions are extended to be global in time for all $p>1$. In the singular case $\dav=0$ and $1<p\le 5$, the mass conservation law guarantees global existence in $L^2(\R)$, and then by the persistence of regularity, $H^1(\R)$-solutions exist globally in time. On the other hand, if the average dispersion $\dav$ is positive, the model admits a much richer dynamics in that it corresponds to the focusing case in the classical NLS. In particular, the model has soliton-like pulses made of ground states.
See \cite{CHLT, CHL, GH, ZGJT01} for the existence and their properties.
They are relatively well explored even in the case $\dav=0$, see \cite{EHL, HKS, HL2009, HL, Kunze, Stanislavova}.

In this paper, we focus on the situation where finite time blowup may occur. We consider the positive average dispersion case $\dav>0$. For numerical simplification, we reduce to the case $\dav=1$, but we also assume that $p\geq 9$. In spite of lack of scaling invariance, we assert that the equation is mass-critical (resp., mass-supercritical) if $p=9$ (resp., $p>9$). Note that given a solution $u(t,x)$ to \eqref{NLS}, $u_\lambda(t,x)=\lambda^{\frac{4}{p-1}}u(\lambda^2 t,\lambda x)$ solves the equation with a different $r$-averaged nonlinearity,
$$
i\partial_t u_\lambda+\partial_x^2u_\lambda
+\int_0^{1/\lambda^2} e^{-ir\partial_x^2}\big(|e^{ir\partial_x^2}u_\lambda|^{p-1}e^{ir\partial_x^2}u_\lambda\big)dr=0
$$
with the initial data $u_{0;\lambda}(x)=\lambda^{\f{4}{p-1}}u_0 (\lambda x)$ for any $\lambda>0$, and that $M[u_{0;\lambda}]=M[u_0 ]$ if and only if $p=9$. Another evidence is that the equation is globally well-posed in $H^1(\R)$ when $1< p<9$, see \cite{CHL}.

\smallskip

The main purpose of our work is to justify that the case $p >9$ can be classified as mass-supercritical for the dispersion managed NLS \eqref{NLS} by showing that it admits finite time blowups. Even more than that, we provide a precise description on the global versus blowup dichotomy for the equation \eqref{NLS}. We note that this is the first result establishing the dichotomy for the dispersion managed NLS. Indeed, even existence of finite time blowup solutions was not known before.

We recall that for the classical power-type NLS, the global versus blowup dichotomy is stated in terms of a specific extremizer for an appropriate Gagliardo-Nirenberg inequality and that a non-scattering solitary wave is given by the extremizer, see \cite{ Holmer, Weinstein}. Then, global solutions in the dichotomy are shown to scatter, so the dynamics below the extremizer are completely characterized, see \cite{Dodson, DodsonMurphy, Holmer_CMP, Killip}. An interesting question is to obtain the same scattering versus blowup dichotomy picture for the dispersion managed NLS. In this direction, our main result gives an affirmative answer for the first step to the so-called Kenig-Merle program.

For the dispersion managed NLS, at first glance, one may guess from the definition of the energy functional that the following {\it local-in-time} Gagliardo-Nirenberg-Strichartz estimate
\beq\label{eq:G-Ntype-0}
\int_0^1\|e^{ir\partial_x^2}f \|_{L^{p+1}}^{p+1}dr \leq C\|f\|_{L^2}^{\f{p+7}{2}}\|\partial_xf\|_{L^2}^{\f{p-5}{2}}
\eeq
would be a candidate substituting the Gagliardo-Nirenberg inequality for the classical NLS. However, due to the presence of the finite interval $[0,1]$, the inequality \eqref{eq:G-Ntype-0} is no longer invariant under scaling. It turns out that the sharp constant for \eqref{eq:G-Ntype-0} is given by that of the {\it global-in-time} Gagliardo-Nirenberg-Strichartz estimate
\beq\label{eq:G-Ntype}
\|e^{ir\partial_x^2}f \|_{L_{r,x}^{p+1}}^{p+1}\leq C\|f\|_{L^2}^{\f{p+7}{2}}\|\partial_xf\|_{L^2}^{\f{p-5}{2}},
\eeq
where
\[
\|u\|_{L^q_{r,x}} := \|u\|_{L_r^q (\R, L_x^q(\R))} = \left( \int_\R \int_\R |u(r,x)|^q dx dr\right)^{\f{1}{q}},\quad 1\leq q<\infty.
\]

\begin{theorem}[Critical element for the inequality \eqref{eq:G-Ntype}]\label{thm:maximizer}
We define the Weinstein functional associated with the inequality \eqref{eq:G-Ntype} by
\beq\label{Weinstein}
W _p(f)=\f{\|e^{ir\partial_x^2}f \|_{L_{r,x}^{p+1}}^{p+1}}{ \|f\|_{L^2}^{\f{p+7}{2}} \|\partial_xf\|_{L^2}^{\f{p-5}{2}}}. \notag
\eeq
Then, for $p\geq 9$, the variational problem
\beq\label{variational on R}
\calC_p=\sup_{f\in H^1(\R)}  W _p(f) \notag
\eeq
admits a maximizer $Q\in H^1(\R)$ which solves the Euler-Lagrange equation
\beq\label{eq:EL}
    -\partial_x^2Q+ Q- \int_\R e^{-ir\partial_x^2}\big(|e^{ir\partial_x^2}Q|^{p-1}e^{ir\partial_x^2}Q\big)dr=0.
\eeq
Moreover, the sharp constant $\calC_p$ for the global inequality \eqref{eq:G-Ntype} is equal to that for the local inequality \eqref{eq:G-Ntype-0}.
\end{theorem}

\begin{remarks}\label{remark: norm quantities}
\begin{theoremlist}
\item Unlike the classical NLS, we do not know uniqueness (up to symmetries) of the critical element for the Weinstein functional. However, by construction, the important norm quantities can be expressed only in terms of $\calC_p$ and $p$, and they are independent of a possibly non-unique profile $Q$. Precisely, we have $\|Q\|_{L^2}^2=\big[\f{2(p+1)(p-5)}{(p+7)^2\calC_p}\big]^{\f{2}{p-1}}\cdot\big(\f{p+7}{p-5}\big)^{\f{1}{2}}$, $\|\partial_xQ\|^2_{L^2}=\f{p-5}{p+7}\|Q\|^2_{L^2}$ and $
\|e^{ir\partial_x^2}Q \|_{L_{r,x}^{p+1}}^{p+1}=\f{2(p+1)}{p+7}\|Q\|_{L^2}^2$, see Lemma \ref{norm identities}.
\item Some precise upper and lower bounds for the sharp constant $\calC_p$ are given in Lemma \ref{upper-lower bound}.
\item As a consequence of Theorem \ref{thm:maximizer}, it is easy to see that there is no maximizer for the local-in-time Gagliardo-Nirenberg-Strichartz estimate \eqref{eq:G-Ntype-0}, see Remark \ref{rem:nonexistence}.
\end{theoremlist}
\end{remarks}

We note that $u(t,x)=Q(x)e^{it}$ is a solitary wave for a rather theoretical scaling-invariant limit equation
\begin{equation}\label{entire r averaged NLS}
i\partial_t u+\partial_x^2u
+\int_{\mathbb{R}} e^{-ir\partial_x^2}\big(|e^{ir\partial_x^2}u|^{p-1}e^{ir\partial_x^2}u\big)dr=0, \notag
\end{equation}
and it preserves the energy over $\R$ defined by
\begin{equation}\label{eq:def_energyQ}
E_\infty[f]:=\f{1}{2}\|\partial_xf\|_{L^2}^2-\f{1}{p+1}\|e^{ir\partial_x^2}f \|_{L_{r,x}^{p+1}}^{p+1} \notag
\end{equation}
and the same mass.

Using the extremizer in the previous theorem, we state our main result on the global existence versus finite time blowup dichotomy.

\begin{theorem}[Global versus blowup criteria]\label{thm:globalvsblowup}
Let $p>9$. Suppose that
\beq\label{ass:1}
E[u_0 ] M[u_0 ]^{\alpha} < E_\infty[Q] M[Q]^{\alpha},
\eeq
where $\alpha=\f{p+7}{p-9}$ and $Q$ is given in Theorem \ref{thm:maximizer}, and let $u(t)$ be the $H^1(\R)$-solution to the equation \eqref{NLS} with initial data $u_0$  whose the maximal time interval of existence is $(-T_{\min},T_{\max})$.

\begin{theoremlist}
\item If \beq\label{ass:global}
\|\partial_xu_0\|_{L^2}\|u_0  \|_{L^2} ^{\alpha}< \|\partial_xQ\|_{L^2}\|Q\|_{L^2}^{\alpha},
\eeq
then 
\beq\label{conc:global}
\|\partial_x u(t)\|_{L^2}\|u(t)\|_{L^2}^{\alpha}< \|\partial_xQ\|_{L^2}\|Q\|_{L^2}^{\alpha} \quad\textup{for all }t\in (-T_{\min},T_{\max})
\eeq
and $u(t)$ exists globally in time, i.e., $T_{\max}=T_{\min}=\infty$.
\item If
\beq\label{ass:blowup}
\|\partial_xu_0  \|_{L^2}\|u_0  \|_{L^2}^{\alpha} > \|\partial_xQ\|_{L^2}\|Q\|_{L^2}^{\alpha}
\eeq
then
\beq\label{conc:blowup}
\|\partial_x u(t)\|_{L^2}\|u(t)\|_{L^2}^{\alpha} >\|\partial_xQ\|_{L^2}\|Q\|_{L^2}^{\alpha}\quad\textup{for all }t\in (-T_{\min},T_{\max}).
\eeq
Furthermore, if $|\cdot|u_0  \in L^2(\R)$, then $T_{\min}<\infty$ so that 
 the solution blows up in finite time in the negative direction.
\end{theoremlist}
\end{theorem}

\begin{remark}
 As mentioned in the introduction, for convenience, the roles of the time and the position variables in our model \eqref{NLS} are switched . The $t$-variable represents the position on the optical cable. Thus, Theorem \ref{thm:globalvsblowup} (ii) indeed establishes a blow-up to the left from the origin.
\end{remark}

\begin{remark}
It follows from Remark \ref{remark: norm quantities} (i) that the key quantities $\|\partial_xQ\|_{L^2}\|Q\|_{L^2}^{\alpha}$ and $E_\infty[Q]M[Q]^{\alpha}$ also can be written using $\calC_p$ and $p$, precisely,
\beq\label{eq: gradient times mass}
\|\partial_xQ\|_{L^2}\|Q\|_{L^2}^{\alpha}= \left[ \f{2(p+1)}{(p-5)\calC_p}\right]^{\f{2}{p-9}}
\eeq
and
\beq\label{eq:energyQ_invariant}
 E_\infty[Q]M[Q]^{\alpha}=\f{p-9}{2(p-5)}\|\partial_xQ\|_{L^2}^{2}\|Q\|_{L^2}^{2\alpha}
 =\f{p-9}{2(p-5)} \left[ \f{2(p+1)}{(p-5)\calC_p} \right]^{\f{4}{p-9}}.
\eeq
\end{remark}

\begin{figure}[htb]\label{figure}
  \centering
  \includegraphics[width=0.65\textwidth]{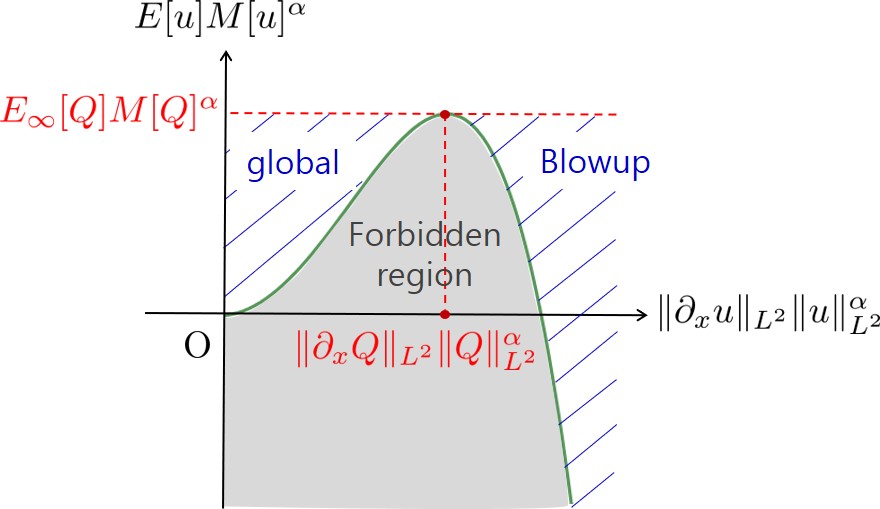}
  \caption{Global existence versus blowup dichotomy}
\end{figure}

\begin{remark}[Negative energy implies blowup]
The subset of $H^1(\mathbb{R})$ satisfying \eqref{ass:1} can be characterized as in Figure \ref{figure}. From the picture, one can see that Theorem \ref{thm:globalvsblowup} (ii) includes that negative energy solutions having finite variance blow up in finite time, which is analogous to the classical result of Glassey \cite{Glassey}.
\end{remark}

\begin{remark}[Mass-critical case, $p=9$]
 In the mass-critical case, global existence is shown in \cite[Theorem 1.6 (ii)]{CHL} when $\|u_0\|_{L^2}$ is small enough. We have a precise sufficient condition
    for global existence in terms of the extremizer in Theorem \ref{thm:maximizer}. More precisely, if $\|u_0 \|_{L^2}< \|Q\|_{L^2}$, then the solution exists global in time. Indeed, it follows from the conservation laws, Theorem \ref{thm:maximizer} and the identities in Remark \ref{remark: norm quantities} (i) that
    \[
    \begin{aligned}
    E[u_0 ]&=E[u(t)]=\f{1}{2}\|\partial_xu(t)\|_{L^2}^2-\f{1}{10} \int_0^1\|e^{ir\partial_x^2}u(t)\|_{L^{10}}^{10} dr\\
    & \geq  \f{1}{2}\|\partial_xu(t)\|_{L^2}^2 -\f{1}{10}\f{\|e^{ir\partial_x^2}Q\|_{L_{r,x}^{10}}^{10}}{ \|Q\|_{L^2}^{8} \|\partial_xQ\|_{L^2}^{2}} \|u(t)\|_{L^2}^{8} \|\partial_xu(t)\|_{L^2}^{2}\\
    &=\f{1}{2}\|\partial_xu(t)\|_{L^2}^2\left(1-\f{\|u_0 \|^8_{L^2}}{\|Q\|_{L^2}^8}\right).
    \end{aligned}
    \]
Therefore, if $\|u_0 \|_{L^2}< \|Q\|_{L^2}$, then $u(t)$ exists globally in time. On the other hand, if $\|u_0 \|_{L^2}>\|Q\|_{L^2}$ and if, in addition, $E[u_0 ]<0$ and $\|xu_0 \|_{L^2}<\infty$, then it follows from Proposition \ref{prop:variance lower bound} and a density argument that the solution blows up in finite time.
\end{remark}

One of the main contributions of this paper is to provide an explicit calculation for the second derivative of the variance for the dispersion managed NLS, see Section \ref{sec:virial}. Similar calculations have potential applications to various dynamical problems of the model in the form of (localized) virial/Morawetz identities, which will be postponed to future work. Indeed, proving virial/Morawetz identities in general can be done by elementary integration by parts with a weight function. However, unlike the power-type case, the dispersion managed NLS includes the linear propagator $e^{ir\partial_x^2}$ in the nonlinearity, so more careful analysis is required to deal with to commute the weight function with the propagator and vice versa, see \eqref{useful identities}. More importantly, while doing integration by parts, additional time boundary terms appear at $r=0, 1$, see \eqref{integralbyparts_r} and \eqref{integralbyparts_r2}. An important remark is that if we repeat the same calculation with the nonlinearity
\[
\int_{r_0}^{r_1} e^{-ir\partial_x^2}(|e^{ir\partial_x^2 } u|^{p-1}e^{ir\partial_x^2 } u ) dr, \quad r_0<r_1,
\]
then a good negative sign term is obtained at the right endpoint $r_1$, while a bad positive sign term is obtained at the left endpoint $r_0$. However, in the physically relevant case, $r_0$ is taken to be zero and it removes the bad sign term. This is a key observation in our proof.

\smallskip

\bigskip

 The paper is organized as follows:
In Section \ref{sec:preliminary}, we review the local theory in $H^1(\R)$ for the Cauchy problem \eqref{NLS} and present the local theory in $H^{3,3}(\R)$ in order to ensure finite variance which is key for the existence of finite time blowup solutions.
 In Section \ref{sec:virial}, we prove the virial estimate by modifying the method of Glassey \cite{Glassey}.
The proof of Theorem \ref{thm:maximizer} which is a central analytic tool in this paper is presented in Section \ref{sec:maximizer}.
Finally, in Section \ref{sec:globalvsblowup}, we establish the global existence versus blowup dichotomy, Theorem \ref{thm:globalvsblowup}.

\bigskip

$\mathbf{Notation.}$
We write $f\lesssim g$ if a finite constant $C>0$ exists such that $f\leq C g$.
 For $k\in \N$, we denote by $H^{k,k}(\R)$ the weighted Sobolev space equipped with the norm
$$
\|f\|_{H^{k,k}}:=\left\{\|f\|_{L^2}^2+\|\partial_x^kf\|_{L^2}^2+\big\|x^kf\big\|_{L^2}^2\right\}^{\frac{1}{2}}.
$$

\section{Preliminary local theory}\label{sec:preliminary}
In this section, we discuss the local theory of \eqref{NLS} in both $H^1(\R)$ and $H^{3,3}(\R)$ for positive times only since the case of negative times is done similarly.
 First, we recall the local existence result in $H^1(\mathbb{R})$ from  \cite{CHL}.

\begin{theorem}[Theorem 1.3 in \cite{CHL}] \label{H1 LWP}
Given initial data $u_0   \in H^1(\R)$, there exists a unique maximal solution
$u\in C([0,T_{\max}); H^1(\R))$ of \eqref{NLS}.  The solution $u(t)$ is maximal in the sense that if $T_{\max}<\infty$ then
$\|u(t)\|_{H^1}\to \infty$ as $t\uparrow T_{\max}$.
Moreover, for all $t\in[0,T_{\max})$, it conserves the mass
\[
M[u(t)]=M[u_0 ]
\]
and the energy
\[
\quad E[u(t)]=E[u_0 ].
\]

\end{theorem}

We complete the local well-posedness in $H^1(\R)$ by proving that solutions depend continuously on initial data. It can be proved by a standard argument, but for the reader's convenience, we give its proof.
\begin{lemma}[Continuity of the data-to-solution map]\label{lemma: continuity of data-to-solution map}
Let $u\in C([0,T_{\textup{max}}); H^1(\mathbb{R}))$ be the maximal solution of \eqref{NLS} with initial data $u_0  \in H^1(\mathbb{R})$. 
Then, for any $T\in (0,T_{\textup{max}})$, there exist small $\delta>0$ and large $K\geq 1$ such that if $\|\tilde{u}_0 -u_0  \|_{H^1}\leq\delta$, then the equation \eqref{NLS} has a unique solution $\tilde{u}\in C([0,T]; H^1(\mathbb{R}))$ with initial data $\tilde{u}_0 \in H^1(\mathbb{R})$ and
$$\|\tilde{u}-u\|_{C([0,T]; H^1)}\leq K\|\tilde{u}_0 -u_0  \|_{H^1}.$$
\end{lemma}

\begin{proof}
It follows from  Proposition 3.5 in \cite{CHL} that there is $T_1>0$, depending only on $p$ and $\|u_0  \|_{H^1}$, such that both $u(t)$ and $\tilde{u}(t)$ exist in $H^1(\R)$ and $\|u(t)\|_{H^1}, \|\tilde u(t)\|_{H^1}\le 2\|u_0  \|_{H^1}$  for all $t \in(0, T_1]$, whenever $\tilde u_0  $ is sufficiently close to $u_0  $ in $H^1(\R)$.

Applying the estimates used in the proof of Theorem 1.3 in \cite{CHL}, one can show that
$$\|(u-\tilde{u})(t)\|_{H^1}\leq \|u_0  -\tilde{u}_0 \|_{H^1} 
+C\int_0^t \Big(\|u(s)\|_{H^1}^{p-1}+\|\tilde{u}(s)\|_{H^1}^{p-1}\Big)\|(u-\tilde{u})(s)\|_{H^1}ds$$
for all $t$.
Hence, it follows from Gr\"onwall's inequality that
\begin{align}
\|(\tilde{u}-u)(t)\|_{H^1}&\leq \|\tilde{u}_0 -u_0  \|_{H^1}\exp\bigg(C\int_0^t \|\tilde{u}(s)\|_{H^1}^{p-1}+\|u(s)\|_{H^1}^{p-1}ds\bigg) \notag\\
&\leq \|\tilde{u}_0 -u_0  \|_{H^1}\exp\bigg(\widetilde{C}\int_0^t \|u(s)\|_{H^1}^{p-1}+\|(\tilde{u}-u)(s)\|_{H^1}^{p-1}ds\bigg). \label{iteration}
\end{align}
Therefore, by taking
$$
K=2\exp\bigg(\widetilde{C}\int_0^T \|u(s)\|_{H^1}^{p-1}ds\bigg)\geq 2,
$$
we can choose $T_*\in (0, T_1]$ such that
\begin{equation}\label{iteration assumption}
\sup_{t\in[0,T_*]}\|(\tilde{u}-u)(t)\|_{H^1}\leq \frac{3K}{4}\|\tilde{u}_0 -u_0  \|_{H^1} \leq \frac{3K\delta}{4} .
\end{equation}

To cover the interval $[0,T]$, we substitute \eqref{iteration assumption} in \eqref{iteration}, then we have
$$\sup_{t\in[0,T_*]}\|(\tilde{u}-u)(t)\|_{H^1}\leq\frac{K}{2}\exp\bigg(\widetilde{C}T\frac{(3K\delta)^{p-1}}{4^{p-1}}\bigg)\|\tilde{u}_0 -u_0  \|_{H^1}.$$
Taking $\delta>0$ small, \eqref{iteration assumption} is improved as
$$\sup_{t\in[0,T_*]}\|(\tilde{u}-u)(t)\|_{H^1}\leq \frac{2K}{3}\|\tilde{u}_0 -u_0  \|_{H^1}.$$
This means that one can extend the existence time for $\tilde{u}(t)$ from $T_*$ to $T$, keeping \eqref{iteration assumption}.
\end{proof}

Next, we prove the local existence in $H^{3,3}(\R)$ employing the following estimate.

\begin{lemma}\label{weighted norm bound}
We have
$$\big\|(x-2it\partial_x)^3 f \big\|_{L^2}\lesssim (1+|t|)^3\|f\|_{H^{3,3}}$$
for all $f\in H^{3,3}(\R)$ and $t \in \R$.
\end{lemma}

\begin{proof}
Since
\[
(x-2it\partial_x)^3=x^3-6itx^2\partial_x-6itx-12t^2x\partial_x^2-12t^2\partial_x+8it^3\partial_x^3,
\]
it suffices to show that
\begin{equation}\label{key inequality for weighted norm bound}
\big\|x^{2}\partial_x f\big\|_{L^2} + \big\|x\partial_x^{2} f\big\|_{L^2}\lesssim \|f\|_{H^{3,3}}
\end{equation}
for all $f\in C_c^\infty(\mathbb{R})$, by density. Taking the integration by parts, we have
$$\begin{aligned}
&\|x^2\partial_x f\|_{L^2}^2+\|x\partial_x^2 f\|_{L^2}^2\\
&=\int_{\mathbb{R}}x^4\partial_x f \overline{\partial_x f}+x^2\partial_x^2 f \overline{\partial_x^2 f}\,dx\\
&=-\int_{\mathbb{R}}\big(4x^3\partial_x f+x^4\partial_x^2 f) \overline{f}\,dx-\int_{\mathbb{R}}(2x\partial_x^2 f+x^2\partial_x^3 f) \overline{\partial_x f}\,dx\\
&=-\int_{\mathbb{R}}\big(4x^3\partial_x f+x^4\partial_x^2 f) \overline{f}\,dx-\int_{\mathbb{R}}x^2\partial_x^3 f \overline{\partial_xf}\,dx+\int_{\mathbb{R}}(2\partial_x^2 f+2x\partial_x^3 f) \overline{f}\, dx.
\end{aligned}$$
Thus, using the Cauchy-Schwarz inequality, we obtain
$$\begin{aligned}
\|x^2\partial_x f\|_{L^2}^2+\|x\partial_x^2 f\|_{L^2}^2&\leq 4\|x^3f\|_{L^2}\|\partial_xf\|_{L^2}+\|x^3f\|_{L^2}\|x\partial_x^2f\|_{L^2}+\|\partial_x^3f\|_{L^2}\|x^2\partial_xf\|_{L^2}\\
&\quad+2\|\partial_x^2f\|_{L^2}\|f\|_{L^2}+2\|\partial_x^3f\|_{L^2}\|xf\|_{L^2}\\
&\leq 8\|f\|_{H^{3,3}}^2+\|f\|_{H^{3,3}}\|x\partial_x^2f\|_{L^2}+\|f\|_{H^{3,3}}\|x^2\partial_xf\|_{L^2}\\
&\leq 9\|f\|_{H^{3,3}}^2+\frac{1}{2}\|x^2\partial_x f\|_{L^2}^2+\frac{1}{2}\|x\partial_x^2 f\|_{L^2}^2
\end{aligned}$$
which yields \eqref{key inequality for weighted norm bound}.
\end{proof}

\begin{lemma}[Local existence in $H^{3,3}(\mathbb{R})$]\label{Hkk LWP}
Given initial data $u_0  \in H^{3,3}(\mathbb{R})$, there exists $\widetilde{T}_{\max}\in (0, \infty]$ and a unique maximal solution $u\in C([0,\widetilde{T}_{\textup{max}}); H^{3,3}(\mathbb{R}))$ of \eqref{NLS}.
The solution $u(t)$ is maximal in the sense that if $\widetilde{T}_{\max}<\infty$, then
$\|u(t)\|_{H^{3,3}}\to \infty$ as $t\uparrow \widetilde{T}_{\max}$.
\end{lemma}
\begin{proof}
For each $T,a >0$,  let
\[
B_{T, a}=\{ u\in C([0, T]; H^{3,3}(\R)) \; :\; \|u\|_{L^\infty([0, T], H^{3,3})}\leq a\}
\]
be equipped with the norm $\|u\|_{L^\infty([0, T]; H^{3,3})}$.
Let $u_0  \in H^{3,3}(\R)$ be fixed and define the map $\Gamma_{u_0  }$ on $B_{T, a}$ by
$$
\Gamma_{u_0 }(u)=e^{it\partial_x^2}u_0  +i\int_0^t\int_0^1 e^{i(t-s)\partial_x^2}e^{-ir\partial_x^2}\big(|e^{ir\partial_x^2}u|^{p-1}e^{ir\partial_x^2}u\big)(s)drds.
$$
Using the identity
\beq\label{identity}
xe^{ir\partial_x}=e^{ir\partial_x} (x-2ir\partial_x), \notag
\eeq
the unitarity of $e^{it\partial_x^2}$ on $L^2(\R)$, and Lemma  \ref{weighted norm bound}, we have
\[
\|\Gamma_{u_0 }(u)\|_{H^{3,3}}\lesssim(1+t)^3\|u_0  \|_{H^{3,3}} +\int_0^t \int_0^1 (1+|t-s-r|)^3 \big\|\big(|e^{ir\partial_x^2}u|^{p-1}e^{ir\partial_x^2}u\big)(s)\big\|_{H^{3,3}}dr ds .
\]
For the nonlinear term, use the fact
\[
x^3|e^{ir\partial_x^2}u|^{p-1}e^{ir\partial_x^2}u=|e^{ir\partial_x^2}u|^{p-1}\Bigl(x^3e^{ir\partial_x^2}u\Bigr)
=|e^{ir\partial_x^2}u|^{p-1}e^{ir\partial_x^2}(x-2ir\partial_x)^3u,
\]
Lemma \ref{weighted norm bound}, the Sobolev embedding $H^1(\R)\hookrightarrow L^\infty(\R)$, and the unitarity of $e^{ir\partial_x^2}$ on $H^1(\R)$, then we obtain
\[
\big\||e^{ir\partial_x^2}u|^{p-1}e^{ir\partial_x^2}u\big\|_{H^{3,3}} \lesssim  \|u\|_{H^1}^{p-1}\|u\|_{H^{3,3}}
\]
for all $r \in[0,1]$. Therefore, combining the two inequalities, we get
\beq \label{eq:H^33}
\|\Gamma_{u_0 }(u)\|_{L^\infty([0,T]; H^{3,3})}\lesssim(1+T)^3 \left( \|u_0  \|_{H^{3,3}} +\int_0^T \|u(s)\|_{H^{3,3}}^{p} ds\right). \notag
\eeq
Similarly, we have
\[
\begin{aligned}
&\|\Gamma_{u_0 }(u)-\Gamma_{u_0 }(v) \|_{L^\infty([0,T]; H^{3,3} )}\\
&\lesssim (1+T)^3\int_0^T\int_0^1\big\|\big(|e^{ir\partial_x^2}u|^{p-1}e^{ir\partial_x^2}u -|e^{ir\partial_x^2}v|^{p-1}e^{ir\partial_x^2}v   \big)(s)\big\|_{H^{3,3}}drds\\
& \lesssim (1+T)^3 \int_0^T \big( \|u(s) \|^{p-1}_{H^{3,3}} +\| v(s)\|^{p-1}_{H^{3,3}}\big)
\|(u-v)(s)\|_{H^{3,3}}ds
\end{aligned}
\]
Thus, for all $u,v \in B_{T,a}$, there exists $C>0$, depending on $p$, such that we have
\[
\|\Gamma_{u_0 }(u)\|_{C([0,T]; H^{3,3})} \leq C (1+T)^3a+C (1+T)^3 Ta^p
\]
and
\[
d(\Gamma_{u_0 }(u), \Gamma_{u_0 }(v)) \leq C (1+T)^3 T a^{p-1}d(u,v).
\]
Taking $a\geq \|u_0  \|_{H^{3,3}}$, it follows that $\Gamma_{u_0 }$ is a contraction from $B_{T,a}$ to itself for a sufficiently small $T>0$, depending on $\|u_0  \|_{H^{3,3}}$.
Therefore, by the Banach fixed point theorem, there exists a unique solution $u$ in $B_{T, a}$.
By a standard argument, we extend the solution $u$ to the maximal interval $[0,\widetilde{T}_{\max})$ so that
$u\in C([0, \widetilde{T}_{\max}); H^{3,3}(\R))$ with either $\widetilde{T}_{\max}=\infty$ or
$\lim_{t\uparrow \widetilde{T}_{\max}} \|u(t)\|_{H^{3,3}} =\infty$.
\end{proof}

We end this section by proving the $H^{3,3}$ regularity of $H^1$ solution.
\begin{lemma}[Persistence of regularity]\label{lemma: persistence of regularity}
Suppose that $u \in C([0,\widetilde{T}_\textup{max}); H^{3,3}(\mathbb{R}))$ is a solution of \eqref{NLS} with initial data $u_0  \in  H^{3,3}(\mathbb{R})$ and the maximal existence time $\widetilde{T}_\textup{max}$ is finite, that is,  $\|u(t)\|_{H^{3,3}}\to\infty$ as $t\uparrow \widetilde{T}_\textup{max}$. Then, $\|u(t)\|_{H^1}\to \infty$ as $t\uparrow \widetilde{T}_\textup{max}$. As a consequence, $\widetilde{T}_{\textup{max}}=T_{\textup{max}}$, where $T_{\textup{max}}$ is the maximal existence time in Theorem \ref{H1 LWP}.
\end{lemma}

\begin{proof}
Suppose to the contrary that $\|u\|_{L^\infty([0,\widetilde{T}_{\textup{max}});H^1)}<\infty$. Arguing similarly as in the proof of Lemma \ref{Hkk LWP}, one can show  that
$$
\|u(t)\|_{H^{3,3}}\lesssim (1+|t|)^3 \|u_0  \|_{H^{3,3}}+\int_0^t (1+|t-s|)^3 \|u(s)\|_{H^1}^{p-1}\|u(s)\|_{H^{3,3}} ds.
$$
Hence, Gr\"onwall's inequality yields the bound
$$
\|u(t)\|_{H^{3,3}}\lesssim (1+|t|)^3 \|u_0  \|_{H^{3,3}}\exp\bigg(\int_0^t (1+|t-s|)^3 \|u(s)\|_{H^1}^{p-1}ds\bigg),
$$
which implies that as long as  $u(t)$ is finite in $H^1(\mathbb{R})$, it is finite in $H^{3,3}(\R)$. Taking the limit as $t\uparrow \widetilde{T}_{\textup{max}}$ deduces a contradiction.
\end{proof}

\section{Virial estimate}\label{sec:virial}

In this section, we modify the standard virial identity of Glassey \cite{Glassey} dealing with the variance. In order to ensure finite variance and use the virial argument, we require the initial data $u_0 $ to be in $H^{3,3}(\R)$.

\begin{proposition}\label{prop:variance lower bound}
Let $u\in C((-\widetilde{T}_{\min}, \widetilde{T}_{\max});H^{3,3}(\R))$ be the maximal solution of \eqref{NLS} with initial data $u_0  \in H^{3,3}(\R)$. If the variance is defined by
$$
\mathcal{V}(t)=\int_{\mathbb{R}}x^2|u(t,x)|^2dx ,
$$
then
\[
\mathcal{V}(t) \le \mathcal{V}(0)+\mathcal{V}_1'(0)t+ \int_0^t\int_0^s \Phi(\tau) d\tau ds
\]
for all $-\widetilde{T}_{\min}<t\leq 0$,
where
\[
\mathcal{V}_1'(t)=4\textup{Im}\int_{\mathbb{R}}x\partial_xu(t,x)\overline{u(t,x)}dx
\]
and
\beq\label{eq:Phi}
\Phi(t)=16E[u_0 ]-\frac{4(p-9)}{p+1}\int_0^1\int_{\mathbb{R}}|e^{ir\partial_x^2}u(t,x)|^{p+1}dxdr.
\eeq
\end{proposition}

\begin{proof}
Let $u\in C((-\widetilde{T}_{\min},\widetilde{T}_{\max});H^{3,3}(\R))$ be the maximal solution of \eqref{NLS}. 
We first calculate the first derivative of $\mathcal{V}$ as follows.
Using equation in \eqref{NLS}, we write
$$\begin{aligned}
\mathcal{V}'(t)&=2\textup{Re}\int_{\mathbb{R}}x^2\overline{u(t,x)}\partial_t u(t,x)dx\\
&=-2\textup{Im}\int_{\mathbb{R}}x^2\overline{u(t,x)}\partial_x^2u(t,x)dx\\
&\quad-2\textup{Im}\int_{\mathbb{R}}\int_0^1 \overline{(e^{ir\partial_x^2}x^2u)(t,x)}\big(|e^{ir\partial_x^2}u|^{p-1}e^{ir\partial_x^2}u\big)(t,x)drdx\\
&=:\mathcal{V}_1'(t)+\mathcal{V}_2'(t).
\end{aligned}$$
By the integration by parts,
$$
\mathcal{V}_1'(t)=4\textup{Im}\int_{\mathbb{R}}x\partial_xu(t,x)\overline{u(t,x)}dx.
$$
Here, the boundary terms disappear since $\partial_x u \in H^1(\R)$ and
$$
\|x^2u\|_{H^1}\lesssim \|xu\|_{L^2}+\|x^2u\|_{L^2}+\|x^2\partial_xu\|_{L^2}\lesssim \|u\|_{H^{3,3}}
$$
where \eqref{key inequality for weighted norm bound} is used.

For $\mathcal{V}_2'(t)$, note the identities
\beq\label{useful identities}
e^{ir\partial_x^2}x=(x+2ir\partial_x)e^{ir\partial_x^2}\quad\textup{and}\quad\partial_x^2 e^{ir\partial_x^2}=-i\partial_r e^{ir\partial_x^2}
\eeq
to get
$$
\begin{aligned}
e^{ir\partial_x^2}x^2&=(x+2ir\partial_x)^2e^{ir\partial_x^2}\\
&=(x^2+4irx\partial_x+2ir-4r^2\partial_x^2)e^{ir\partial_x^2}\\
&=\big(x^2+i(4rx\partial_x+2r+4r^2\partial_r)\big)e^{ir\partial_x^2}.
\end{aligned}
$$
 Hence, using the integration by parts again, we obtain
\beq\label{integralbyparts_r}
\begin{aligned}
\mathcal{V}_2'(t)&=2\textup{Re}\int_{\mathbb{R}}\int_0^1 \overline{\big((4rx\partial_x+2r+4r^2\partial_r)e^{ir\partial_x^2}u\big)(t,x)}(|e^{ir\partial_x^2}u|^{p-1}e^{ir\partial_x^2}u)(t,x)drdx\\
&=2\int_\R \int_0^1 \left(\left(\frac{4}{p+1} rx\partial_x+2r +\frac{4r^2}{p+1}\partial_r\right)|e^{ir\partial_x^2} u|^{p+1}\right) (t,x)drdx \\
&=\frac{4(p-5)}{p+1}\int_{\mathbb{R}}\int_0^1|e^{ir\partial_x^2}u(t,x)|^{p+1}rdrdx  +\frac{8}{p+1}\int_0^1\left[x|e^{ir\partial_x^2}u(t,x)|^{p+1}\Big|_{x=-\infty}^{x=\infty}\right]rdr \\
&\quad +\frac{8}{p+1}\int_{\mathbb{R}}|e^{ir\partial_x^2}u(t,x)|^{p+1}r^2\Big|_{r=0}^{r=1} dx.
\end{aligned}
\eeq
Now we note that $x|e^{ir\partial_x^2}u(t,x)|^{p+1}$ is in $H^1(\R)$ for each $r\in[0,1]$.  Indeed, by the Sobolev inequality,  \eqref{useful identities}, and a similar argument as when obtaining \eqref{key inequality for weighted norm bound}, we see that
\[
\begin{aligned}
\|x|e^{ir\partial_x^2}u|^{p+1}\|_{L^2}&\le \|e^{ir\partial_x^2}u\|^p_{L^\infty}\|xe^{ir\partial_x^2}u\|_{L^2}\\
&\le \|e^{ir\partial_x^2}u\|^p_{H^1}\|(x-2ir\partial_x)u\|_{L^2}
\lesssim (1+r)\|u\|^p_{H^1}\|u\|_{H^{1,1}}
\end{aligned}\]
and
\[
\begin{aligned}
\|x|e^{ir\partial_x^2}u|^{p}\partial_x(e^{ir\partial_x^2}u)\|_{L^2}
& \le \|e^{ir\partial_x^2}u\|^p_{L^\infty}\|xe^{ir\partial_x^2}(\partial_x u)\|_{L^2}\\
&\le \|e^{ir\partial_x^2}u\|^p_{H^1}\|(x-2ir\partial_x)(\partial_x u)\|_{L^2}
\lesssim (1+r)\|u\|^p_{H^1}\|u\|_{H^{2,2}}
\end{aligned}
\]
and therefore $x|e^{ir\partial_x^2}u(t,x)|^{p+1} \in H^1(\R)$ since $u \in H^{3,3}(\R)$. Thus, the boundary terms $x|e^{ir\partial_x^2}u(t,x)|^{p+1}\Big|_{-\infty}^\infty$ vanish and so
\begin{equation}\label{V_2}
\mathcal{V}_2'(t)=\frac{4(p-5)}{p+1}\int_{\mathbb{R}}\int_0^1|e^{ir\partial_x^2}u(t,x)|^{p+1}rdrdx
+\frac{8}{p+1}\int_{\mathbb{R}}|e^{i\partial_x^2}u(t,x)|^{p+1}dx>0.
\end{equation}
Therefore, using the fundamental theorem of calculus, we have
\begin{equation}\label{V_1}
\mathcal{V}(0)-\mathcal{V}(t)= \int_t^0 \mathcal{V}'(\tau)d\tau \ge \int_t^0 \mathcal{V}_1'(\tau)d\tau , \quad\mbox{i.e.,}\quad \mathcal{V}(t)\leq \mathcal{V}(0)-\int_t^0 \mathcal{V}_1'(\tau)d\tau
\end{equation}
for all $t<0$ since $\mathcal{V}'(t)=\mathcal{V}_1'(t)+\mathcal{V}_2'(t) > \mathcal{V}_1'(t)$ by \eqref{V_2}.

Differentiating $\mathcal{V}_1'(t)=4\textup{Im}\int_{\mathbb{R}}x\partial_xu(t,x)\overline{u(t,x)}dx$ and then using the integration by parts, we obtain
$$
\begin{aligned}
\mathcal{V}_1''(t)
&=4\textup{Im}\int_{\mathbb{R}}x \Bigl(\partial_x\partial_t u(t,x)\Bigr)
\overline{u(t,x)}dx+4\textup{Im}\int_{\mathbb{R}}x\partial_xu(t,x)\overline{\partial_tu(t,x)}dx\\
&=-4\textup{Im}\int_{\mathbb{R}}\partial_t u(t,x)\overline{u(t,x)}dx+8\textup{Im}\int_{\mathbb{R}}x\partial_xu(t,x)\overline{\partial_tu(t,x)}dx .
\end{aligned}
$$
Then, using equation in \eqref{NLS} and the integration by parts, again, we get
$$\begin{aligned}
\mathcal{V}_1''(t)&=-4\textup{Re}\int_{\mathbb{R}}\left\{\partial_x^2u+\int_0^1 e^{-ir\partial_x^2}(|e^{ir\partial_x^2}u|^{p-1}e^{ir\partial_x^2}u)dr\right\}(t,x)\overline{u(t,x)}dx\\
&\quad-8\textup{Re}\int_{\mathbb{R}}x\left\{\partial_x^2u+\int_0^1 e^{-ir\partial_x^2}(|e^{ir\partial_x^2}u|^{p-1}e^{ir\partial_x^2}u)dr\right\}(t,x)\overline{\partial_xu(t,x)}dx\\
&=8\|\partial_xu(t,x)\|_{L^2}^2-4\int_{\mathbb{R}}\int_0^1|e^{ir\partial_x^2}u(t,x)|^{p+1}drdx\\
&\quad-8\textup{Re}\int_{\mathbb{R}}\int_0^1(|e^{ir\partial_x^2}u|^{p-1}e^{ir\partial_x^2}u)(t,x)\overline{e^{ir\partial_x^2}(x\partial_xu)(t,x)}drdx.
\end{aligned}$$
On the other hand, it follows from the identities in \eqref{useful identities} that
\beq\label{eq:identity}
e^{ir\partial_x^2}x\partial_x=(x+2ir\partial_x) e^{ir\partial_x^2}\partial_x=x\partial_xe^{ir\partial_x^2}+2r\partial_r e^{ir\partial_x^2}. \notag
\eeq
Hence, by a similar argument in \eqref{integralbyparts_r}, we have
\beq\label{integralbyparts_r2}
\begin{aligned}
&-\textup{Re}\int_{\mathbb{R}}\int_0^1(|e^{ir\partial_x^2}u|^{p-1}e^{ir\partial_x^2}u)(t,x)\overline{e^{ir\partial_x^2}(x\partial_xu)(t,x)}drdx\\
&=-\textup{Re}\int_{\mathbb{R}}\int_0^1(|e^{ir\partial_x^2}u|^{p-1}e^{ir\partial_x^2}u)(t,x)\overline{(x\partial_xe^{ir\partial_x^2}u+2r\partial_re^{ir\partial_x^2}u)(t,x)}drdx\\
&= -\f{1}{p+1}\int_\R\int_0^1 \left((x\partial_x+2r\partial_r)|e^{ir\partial_x^2}u|^{p+1}\right)(t,x)drdx\\
&=\frac{3}{p+1}\int_{\mathbb{R}}\int_0^1|e^{ir\partial_x^2}u(t,x)|^{p+1}drdx-\frac{2}{p+1}\int_{\mathbb{R}}|e^{i\partial_x^2}u(t,x)|^{p+1}dx\\
&\leq \frac{3}{p+1}\int_{\mathbb{R}}\int_0^1|e^{ir\partial_x^2}u(t,x)|^{p+1}drdx.
\end{aligned}
\end{equation}
Therefore, it follows from the energy conservation that
\[ \label{ineq:v''}
\begin{aligned}
\mathcal{V}_1''(t)&\leq8\|\partial_xu(t)\|_{L^2}^2
-\frac{4(p-5)}{p+1}\int_0^1\int_{\mathbb{R}}|e^{ir\partial_x^2}u(t,x)|^{p+1}dxdr\\
&=16E[u_0 ]-\frac{4(p-9)}{p+1}\int_0^1\int_{\mathbb{R}}|e^{ir\partial_x^2}u(t,x)|^{p+1}dxdr =\Phi(t).
\end{aligned}
\]
Thus, by the fundamental theorem of calculus again, for all $t<0$
\[
\mathcal{V}_1'(0)-\mathcal{V}_1'(t)=\int_t^0 \mathcal{V}_1''(s) ds \leq \int_t^0 \Phi(s) ds ,
\quad\mbox{i.e.,}\quad 
\mathcal{V}_1'(t) \ge \mathcal{V}_1'(0)-\int_t^0 \Phi(s) ds
\]
and therefore by \eqref{V_1}
\[
\mathcal{V}(t) \le \mathcal{V}(0)+\mathcal{V}_1'(0)t+\int_t^0 \int_s^0  \Phi(\tau)d\tau ds.
\]

\end{proof}

\section{Proof of Theorem \ref{thm:maximizer}}\label{sec:maximizer}

In this section, we consider the variational problem
\beq\label{eq:variational}
\calC_p=\sup_{f\in H^1(\R)}W_p(f)= \sup_{f\in H^1(\R)} \f{ \|e^{ir\partial_x^2}f \| _{L_{r,x}^{p+1}}^{p+1} }{\|f\|_{L^2}^{\f{p+7}{2}} \|\partial_x f\|_{L^2}^{\f{p-5}{2}}}, \notag
\eeq
where $\|u\|_{L^q_{r,x}} := \|u\|_{L_r^q (\R, L_x^q(\R))}$ for $1\leq q<\infty$. The following lemma provides a rough estimate for the best constant $\calC_p$, so the variational problem $\calC_p$ is well-formulated.

\begin{lemma}[Upper and lower bounds for $\calC_p$]\label{upper-lower bound}
\[
\f{2^{\f{p-7}{4}}}{\pi^{\f{p-1}{4}}(p+1)^{\f{1}{2}}}  B\left(\f{1}{2}, \f{p-3}{4}\right) = \f{2^{\f{p-7}{4}}} {\pi^{\f{p-3}{4}}(p+1)^{\f{1}{2}}} \f{\Gamma(\f{p-3}{4})}{\Gamma(\f{p-1}{4})}\leq \calC_p\leq \frac{1}{2\sqrt{3}},
\]
where $B$ is the beta function and $\Gamma$ is the gamma function. In particular, $\frac{1}{2\sqrt{5}\pi} \leq  \calC_9\leq \frac{1}{2\sqrt{3}}$.
\end{lemma}

\begin{proof}
For every $f\in H^1(\R)$, we have
\[ \label{ineq:G-NandStraichartz}
\begin{aligned}
 \|e^{ir\partial_x^2}f \| _{L_{r,x}^{p+1}}^{p+1}&\leq
 \int_\R \|e^{ir\partial_x^2}f \|_{L^\infty}^{p-5} \|e^{ir\partial_x^2}f \|_{L^6}^{6} dr \\
& \leq \|f\|_{L^2}^{\f{p-5}{2}}\|\partial_x f\|_{L^2}^{\f{p-5}{2}}\|e^{ir\partial_x^2}f \|_{L_{r,x}^6}^6.
\end{aligned}
\]
Then, by the Strichartz estimate in one dimension
\[
\|e^{ir\partial_x^2}f \|_{L_{r,x}^6} \leq \frac{1}{12^{1/12}} \|f\|_{L^2}
\]
with the best constant $\frac{1}{12^{1/12}}$, see \cite{Foschi, HZ}, it follows that
\beq\label{ineq:G-N_type_R}
\|e^{ir\partial_x^2}f \| _{L_{r,x}^{p+1}}^{p+1} \leq
\frac{1}{\sqrt{12}}\|f\|_{L^2}^{\f{p+7}{2}} \|\partial_x f\|_{L^2}^{\f{p-5}{2}}.
\eeq

For the lower bound, we employ the Gaussian test-function $\varphi(x)=e^{-\f{1}{2}x^2}$. Note that
$\|\varphi\|_{L^2}^2=\sqrt{\pi}$ and $\|\partial_x\varphi\|_{L^2}^2=\sqrt{\pi}/2$.
Moreover, its free time evolution is given
by
\[
e^{ir\partial_x^2}\varphi(x)=\left( \f{1}{1+2ir}\right)^{\f{1}{2}}e^{-\f{x^2}{2(1+2ir)}}
\]
and therefore
\[
\|e^{ir\partial_x^2}\varphi\|_{L^{p+1}}^{p+1}=\left(\f{2\pi}{p+1}\right)^{\f 1 2}\left(\f{1}{1+4r^2}\right)^{\f {p-1}{4}}
\]
see, e.g., \cite{CHLT, ZGJT01}. Then we have
\[
W_p(\varphi)=\f{\|e^{ir\partial_x^2}\varphi\|_{L_{r,x}^{p+1}}^{p+1}}{\|\varphi\|_{L^2}^{\f{p+7}{2}}\|\partial_x\varphi\|_{L^2}^{\f{p-5}{2}}}
=\f{2^{\f{p+1}{4}}} {\pi^{\f{p-1}{4}}(p+1)^{\f{1}{2}}} \int_0^\infty \left(\f{1}{1+4r^2}\right)^{\f{p-1}{4}} dr.
\]
The integral is rewritten as the beta function that
\[
\int_0^\infty \left(\f{1}{1+4r^2}\right)^{\f{p-1}{4}} dr= \f{1}{2}\int_0^{\f{\pi}{2}} \cos^{\f{p-5}{2}}\theta d\theta =\f{1}{4}B\left(\f{1}{2}, \f{p-3}{4}\right)
\]
where we used $B(x,y)=2\int_0^{\pi/2} (\sin\theta)^{2x-1} (\cos\theta)^{2y-1} d\theta$.
Using the relation between the beta function and gamma function, we obtain the lower bound in the lemma.
\end{proof}

We prove the existence of a maximizer for the variational problem $\calC_p$ and derive its Euler-Lagrange equation.

\begin{proposition}[Existence of a maximizer for $\calC_p$ ]\label{prop:maximizer}
For $p\geq 9$, there exists a maximizer $ g \in H^1(\R)$ for  the variational problem $\calC_p$ such that $\|  g \| _{L^2}=\| \partial_x g\|_{L^2}=1$.
Moreover, $ g $ solves
\begin{equation}\label{eq:EL-tilde}
\f{p-5}{2} \partial_x^2 g-\f{p+7}{2} g
+\f{p+1}{\calC_p} \int_\R e^{-ir\partial_x^2}\big(|e^{ir\partial_x^2} g |^{p-1}e^{ir\partial_x^2} g \big)dr=0. \notag
\end{equation}
\end{proposition}

To begin with, such a critical element can be constructed by a standard concentration-compactness argument of Lions \cite{Lions1, Lions2}. Among many different formulations, we particularly employ one in the form of the profile decomposition, see \cite[Proposition 3.4]{Guevara} for instance.

\begin{lemma}[Profile decomposition]\label{lem:decomposition}
Let $\{f_n\}_{n\in \N}$ be a uniformly bounded sequence in $H^1(\R)$ such that $\displaystyle\limsup_{n\to\infty} \|f_n\|_{L^2} \neq 0$. Then, up to a subsequence, there exists a nonzero $ g  \in H^1(\R)$ and sequences of time shifts $\{r_n\}_{n\in \N}$ and space shifts $\{x_n\}_{n\in \N}$ such that
\[
e^{ir_n \partial_x^2} f_n(\cdot +x_n )\rightharpoonup  g   \quad\mbox{in}\quad H^1(\R).
\]
We define the sequence $R_n$ of remainders by
\[
f_n=  e^{-ir_n \partial_x^2} g (\cdot-x_n)+R_n.
\]
Then, the above decomposition satisfies the asymptotic Pythagorean rule,
\beq\label{eq:profile-1}
\|f_n\|_{L^2}^2=\| g \|_{L^2}^2+ \|R_n\|_{L^2}^2+o_n(1),
\eeq
\beq\label{eq:profile-2}
\|\partial_x f_n\|_{L^2}^2=\| \partial_x g\|_{L^2}^2+ \|\partial_x R_n\|_{L^2}^2+o_n(1)
\eeq
and
\beq\label{eq:profile-3}
\|e^{ir\partial_x^2}f_n\|_{L_{r,x}^{p+1}}^{p+1}=\|e^{ir\partial_x^2} g \|_{L_{r,x}^{p+1}}^{p+1}+ \|e^{ir\partial_x^2}R_n\|_{L_{r,x}^{p+1}}^{p+1}+o_n(1).
\eeq
\end{lemma}
We need to show \eqref{eq:profile-3} only, as all the other can be obtained from \cite[Propositions 3.4]{Guevara} and its proof.
We modify the method by Br\'{e}zis-Lieb \cite{BL} to prove \eqref{eq:profile-3}.

\begin{proof}
Replacing $e^{ir_n \partial_x^2} f_n(\cdot +x_n )$ by $f_n$, we may assume that $x_n=0$ and $r_n=0$ so that $f_n\rightharpoonup  g $ in $H^1(\R)$ and $f_n= g +R_n$, where $R_n=e^{ir_n \partial_x^2} R_n(\cdot +x_n )$ . Then, there exists a subsequence of $\{f_n\}_{n\in\N}$, but still denoted by $\{f_n\}_{n\in\N}$,  such that $f_n \to g $ a.e. in $\R$.
It suffices to prove that
   \beq\label{conv:splitting}
    \lim_{n\to \infty} \int_\R\int_{\R} \left| |e^{ir\partial_x^2}f_{n}|^{p+1} -|e^{ir\partial_x^2}R_n|^{p+1} - |e^{ir\partial_x^2}  g |^{p+1}\right| dxdr =0.
    \eeq
Let $\veps>0$ be fixed. Observe that by Young's inequality
    \beq\label{eq:Young's}
    \begin{aligned}
    ||a+b|^{p+1}-|a|^{p+1}|&=\left|\int_0^1 \f{d}{ds} (|a+sb|^{p+1})ds\right|\\
    &\leq ({p+1}) \int_0 ^1 |a+ sb|^{p} |b| ds\\
    & \leq C\veps |a|^{p+1} +C(\veps) |b|^{p+1}
      \end{aligned}
    \eeq
    for all $a, b \in \C$.
Applying \eqref{eq:Young's} with $a=e^{ir\partial_x^2}R_n$ and $b=e^{ir\partial_x^2} g $, we see that
 \beq\label{ineq:aandb}
 \begin{aligned}
&\left| |e^{ir\partial_x^2}f_n|^{p+1} -|e^{ir\partial_x^2}R_n|^{p+1} - |e^{ir\partial_x^2}  g |^{p+1}\right| \\
& \leq \left| |e^{ir\partial_x^2}f_n|^{p+1} -|e^{ir\partial_x^2}R_n|^{p+1}\right| +|e^{ir\partial_x^2}  g |^{p+1}\\
 & \leq C \veps |e^{ir\partial_x^2} R_n|^{p+1} +C(\veps)|e^{ir\partial_x^2}  g |^{p+1} + |e^{ir\partial_x^2} g |^{p+1}.
 \end{aligned}
 \eeq
Let
\[
F_{n}^\veps(r,x):= \left[ \left| |e^{ir\partial_x^2}f_n(x)|^{p+1} -|e^{ir\partial_x^2}R_n(x)|^{p+1} - |e^{ir\partial_x^2}  g (x)|^{p+1}\right| -C \veps |e^{ir\partial_x^2} R_n(x)|^{p+1}\right]_+,
\]
where $a_+=\max(a,0)$ for any $a\in\R$. For now, we fix $r\in\R$. Then, up to a subsequence,
$F_{n}^\veps (r,\cdot)\to 0$  for a.e. $x$ since $\{e^{ir\partial_x^2}f_n\}_{n\in\N}$ is also uniformly bounded in $H^1(\R)$ and $f_n \to  g $ for a.e. $x$.
Moreover, by \eqref{ineq:aandb} and the Sobolev embedding,
\[
F_{n}^\veps(r,\cdot)\leq C(\veps)|e^{ir\partial_x^2}  g |^{p+1} + |e^{ir\partial_x^2} g |^{p+1}  \in L_{x}^1(\R).
\]
Thus, by the dominated convergence theorem,
$\int_\R F_{n}^\veps(r,x) \, dx \to 0$.
Next, the sequence of functions of $r$ is bounded as
\[
\int_\R F_{n}^\veps(r,x) \, dx \leq \int_\R \left( C(\veps)|e^{ir\partial_x^2}  g |^{p+1} + |e^{ir\partial_x^2} g |^{p+1} \right) dx  \in L_{r}^1(\R)
\]
by \eqref{ineq:G-N_type_R}.
Using the dominated convergence theorem again,  we have $\int_\R\int_\R F_{n}^\veps(r,x) \, dxdr \to 0$.
Hence, we have
 \[
 \begin{aligned}
&\limsup_{n\to \infty}\int_\R\int_\R  \left| |e^{ir\partial_x^2}f_n|^{p+1} -|e^{ir\partial_x^2}R_n|^{p+1} - |e^{ir\partial_x^2}  g |^{p+1}\right|dxdr\\
 &\leq \limsup_{n\to \infty} \left(  \int_\R \int_\R  F_{n}^\veps(r) dxdr  + C\veps \int_\R \int_\R|e^{ir\partial_x^2} R_n|^{p+1} dxdr \right) \leq C\veps
 \end{aligned}
\]
by $\{R_n\}_{n\in\N}$ begin uniformly bounded in $H^1(\R)$ and \eqref{ineq:G-N_type_R}.
Since $\veps>0$ is arbitrary, this proves \eqref{conv:splitting}.
\end{proof}

We use the following elementary inequality to eliminate the splitting scenario in the concentration-compactness argument.
\begin{lemma}[Algebraic inequality \cite{Hong}]\label{lem:algebraic}
  For any $a_1, a_2, b_1, b_2>0$ and $s>1$, if $\f{a_1}{a_2^s}\geq \f{b_1}{b_2^s}$, then
\beq\label{algebraic inequality}
\f{a_1}{a_2^s} \geq
\f{a_1+b_1 }{(a_2+b_2)^s}. \notag
\eeq
Further assume that $\delta\leq \f{b_2}{a_2}\leq \f{1}{\delta}$ for some $\delta>0$, then there exists
$c=c(\delta)>1$ such that
\[
\f{a_1}{a_2^s} \geq
c\f{a_1+b_1 }{(a_2+b_2)^s}.
\]
\end{lemma}

\bigskip

Now we are ready to give the
\begin{proof}[Proof of Proposition  \ref{prop:maximizer}]
Let $\{f_n\}_{n\in \N}\subset H^1(\R)$ be a maximizing sequence for $\calC_p$, i.e., $W_p(f_n)\to\calC_p$.
Note that if define
\begin{equation}\label{eq:f_lambda,mu}
f_{\lambda, \mu}(x)=\mu f(\lambda x) \quad\mbox{for any }\lambda, \mu >0,
\end{equation}
then
\beq\label{eq:scaling}
\|f_{\lambda, \mu}\|_{L^2}^2=\mu^2 \lambda^{-1}\|f\|^2_{L^2},\quad
\|\partial_xf_{\lambda, \mu}\|^2_{L^2}=\mu^2 \lambda\|\partial_xf\|^2_{L^2}
\eeq
and
\[
\|e^{ir\partial_x^2}f_{\lambda, \mu} \| _{L^{p+1}_{r,x}}^{p+1}=\mu^{p+1}\lambda^{-3} \|e^{ir\partial_x^2}f\| _{L^{p+1}_{r,x}}^{p+1}
\]
and therefore  $W_\R^p(f_{\lambda, \mu})=W_\R^p(f)$, that is, the functional $W_{\R}^p(f)$ is invariant under scaling and multiplication by a constant. Thus, we may assume that $\|f_n\|_{L^2}=\|f_n'\|_{L^2}=1$ for each $n\in\N$.

According to Lemma \ref{lem:decomposition},
there is a subsequence of $\{f_n\}_{n\in \N}$, but still denoted by $\{f_n\}_{n\in \N}$, a nonzero profile $ g \in H^1(\R)$, sequences of time shifts $\{r_n\}_{n\in \N}$ and space shifts $\{x_n\}_{n\in \N}$ such that
\beq\label{eq: profile decomposition}
f_n= e^{-ir_n \partial_x^2} g (\cdot-x_n)+ R_n
\eeq
satisfying \eqref{eq:profile-1}-\eqref{eq:profile-3}. Note from \eqref{eq:profile-1} and \eqref{eq:profile-2} that $0<\| g \|_{L^2},\| \partial_x g\|_{L^2}\leq 1$.
Moreover, by the space-time translation invariance of the quantities in the functional $W_p(f)$, we may assume that  $r_n=0$ and $x_n=0$ for all $n$ in \eqref{eq: profile decomposition} so that
\begin{equation}\label{eq: profile decomposition'}
f_n= g + R_n. \notag
\end{equation}

We claim that $ g $ is a maximizer with $\| g \|_{L^2}=\| \partial_x g\|_{L^2}=1$. To prove the claim, it suffices to show
\begin{equation}\label{eq:remainder}
\liminf_{n\to\infty}\|e^{ir\partial_x^2}R_n \|_{L^{p+1}_{r,x}}= 0 .
\end{equation}
Indeed, it follows from \eqref{eq:profile-3} with \eqref{eq:remainder} and $0<\| g \|_{L^2},\| \partial_x g\|_{L^2}\leq 1$ that
\[
\begin{aligned}
\calC_p&=\lim_{n\to\infty} W_p(f_n) =\lim_{n\to\infty} \|e^{ir\partial_x^2}f_n\| _{L_{r,x}^{p+1}}^{p+1}=\|e^{ir\partial_x^2} g \| _{L_{r,x}^{p+1}}^{p+1}\\
&\leq \f{ \|e^{ir\partial_x^2} g \| _{L_{r,x}^{p+1}}^{p+1}}{\| g \|_{L^2}^{\f{p+7}{2}} \| \partial_x g\|_{L^2}^{\f{p-5}{2}}}= W_p( g )\leq \calC_p .
\end{aligned}\]
Therefore, we conclude that $W_p( g )=\calC_p $ and $\| g \|_{L^2}=\| \partial_x g\|_{L^2}=1$. Now, for \eqref{eq:remainder}, suppose to the contrary that
\[\label{eq: dichotomy scenario}
\liminf_{n\to\infty}\|e^{ir\partial_x^2}R_n \|_{L^{p+1}_{r,x}}>0.
\]
Then, it follows from \eqref{ineq:G-N_type_R}, \eqref{eq:profile-1}, and \eqref{eq:profile-2} that
\begin{equation}\label{eq: dichotomy scenario'}
0< \lim_{n\to\infty}\|R_n\|_{L^2}\leq 1,\quad 0< \lim_{n\to\infty}\|\partial_x R_n\|_{L^2}\leq 1.
\end{equation}
Using \eqref{eq:profile-1}-\eqref{eq:profile-3}, we can write
\[\label{eq: W(fn) decomposition}
\begin{aligned}
W_p(f_n)=\f
{ \|e^{ir\partial_x^2} g \|_{L_{r,x}^{p+1}}^{p+1}  + \| e^{ir\partial_x^2} R_n \| _{L_{r,x}^{p+1}}^{p+1} }{ \left(\| g \|_{L^2}^2+\|R_n\|_{L^2}^2\right)^{\f{p+7}{4}} } +o_n(1).
\end{aligned}
\]
We consider the following two cases,
\[
\f{ \|e^{ir\partial_x^2} g  \| _{L_{r,x}^{p+1}}^{p+1} }{\| g \|_{L^2}^{\f{p+7}{2}}} \geq \liminf_{n\to\infty}\f{ \|e^{ir\partial_x^2}R_n \| _{L_{r,x}^{p+1}}^{p+1} }{ \|R_n\|_{L^2}^{\f{p+7}{2}}}
\quad \mbox{and}\quad
\f{ \|e^{ir\partial_x^2} g  \| _{L_{r,x}^{p+1}}^{p+1} }{\| g \|_{L^2}^{\f{p+7}{2}}} < \liminf_{n\to\infty}\f{ \|e^{ir\partial_x^2}R_n \| _{L_{r,x}^{p+1}}^{p+1} }{ \|R_n\|_{L^2}^{\f{p+7}{2}}}.
\]
In the former case, by Lemma \ref{lem:algebraic}, we have
\[
W_p(f_n)\leq \f{ \|e^{ir\partial_x^2} g  \| _{L_{r,x}^{p+1}}^{p+1} }{\| g \|_{L^2}^{\f{p+7}{2}}} +o_n(1)
\]
up to a subsequence. Then, since $0<\| \partial_x g\|_{L^2}\leq 1$, we obtain
\[
\calC_p=\lim_{n\to \infty}W_p(f_n)\leq \f{ \|e^{ir\partial_x^2} g  \| _{L_{r,x}^{p+1}}^{p+1} }{\| g \|_{L^2}^{\f{p+7}{2}}}\leq \f{ \|e^{ir\partial_x^2} g  \| _{L_{r,x}^{p+1}}^{p+1} }{\| g \|_{L^2}^{\f{p+7}{2}} \| \partial_x g\|_{L^2}^{\f{p-5}{2}}} =W_p( g )\leq \calC_p,
\]
 which implies $\| \partial_x g\|_{L^2}=1$. It follows from \eqref{eq:profile-2} and $\|\partial_x f_n\|_{L^2}=1$ that $\|\partial_x R_n\|_{L^2}\to 0$, which contradicts \eqref{eq: dichotomy scenario'}.
In the latter case, let $\delta=\lim_{n\to\infty}\|R_n\|_{L^2}^2$ then $0<\delta \le 1$  due to \eqref{eq: dichotomy scenario'}.
By applying Lemma \ref{lem:algebraic}, there exists $c>1$, independent of $n$, such that
\[
c W_p(f_n) \leq \, \f{ \|e^{ir\partial_x^2}R_n \| _{L_{r,x}^{p+1}}^{p+1} }{\|R_n\|_{L^2}^{\f{p+7}{2}}}+o_n(1) \leq  W_p(R_n) +o_n(1)
\]
for sufficiently large $n$, where we used \eqref{eq: dichotomy scenario'} in the second inequality.
Thus, one sees that
\[
\calC_p =\lim_{n\to \infty } W_p(f_n) < c\lim_{n\to \infty } W_p(f_n) \le \liminf_{n\to \infty } W_p(R_n),
\]
which is a contradiction.

From the standard argument in the calculus of variations, the above maximizer $ g \in H^1(\R)$ is a weak solution of the associated Euler-Lagrange equation
\[
\left.\f{d}{d\varepsilon}\right|_{\varepsilon=0}W_p( g +\varepsilon \eta)=0
\]
for all $\eta\in C_0^\infty(\R)$.
Considering that $\| g \|_{L^2}=1$ and $\| \partial_x g\|_{L^2}=1$, we have
\beq\label{eq:EL_g}
\f{p-5}{2} \partial_x^2 g-\f{p+7}{2} g
+\f{p+1}{\calC_p} \int_\R e^{-ir\partial_x^2}\big(|e^{ir\partial_x^2} g |^{p-1}e^{ir\partial_x^2} g \big)dr=0.
\eeq
Observe that every weak solution $f\in H^1(\R)$ of \eqref{eq:EL_g} is in $H^3(\R)$, since
\[
f \mapsto \int_\R e^{-ir\partial_x^2}\big(|e^{ir\partial_x^2} f |^{p-1}e^{ir\partial_x^2} f \big)dr
\]
maps $H^1(\R)$ into itself. 
Thus, the maximizer $g\in H^3(\R)$ is a strong solution for \eqref{eq:EL_g}.
\end{proof}

From the maximizer we constructed in Proposition \ref{prop:maximizer}, we find a maximizer for $\calC_p$ whose Euler-Lagrange equation is given in \eqref{eq:EL} and we prove that $\calC_p$ equals the sharp constant for the local-in-time Gagliardo-Nirenberg-Strichartz estimate.
\begin{proof}[Proof of Theorem \ref{thm:maximizer}]

We define $Q(x)=\mu  g (\lambda x)$, where $g$ is the maximizer for the variational problem $\calC_p$ constructed in Proposition \ref{prop:maximizer},
\begin{equation}\label{eq: lambda and mu choice}
\lambda=\left(\f{p-5}{p+7}\right)^{\f{1}{2}}\quad\mbox{and}\quad
\mu=\left[\f{2(p+1)(p-5)}{(p+7)^2\calC_p}\right]^{\f{1}{p-1}}. \notag
\end{equation}
Note that $Q$ is also a maximizer, because the functional $W_p(f)$ is invariant under scaling and multiplication by a constant. Moreover, by direct calculations, one can see that this modified critical element $Q$ solves the equation with normalized coefficients
\beq\label{eq:scalingE-L}
 -\partial_x^2Q+ Q- \int_\R e^{-ir\partial_x^2}\big(|e^{ir\partial_x^2}Q|^{p-1}e^{ir\partial_x^2}Q\big)dr=0 \notag
\eeq
since
\begin{equation}\label{eq:scaling-nonlinearity}
\int_\R e^{-ir\partial_x^2}\big(|e^{ir\partial_x^2}\mu  g (\lambda x)|^{p-1}e^{ir\partial_x^2}\mu  g (\lambda x)\big)dr
= \mu^p \lambda^{-2}\int_\R e^{-ir\partial_x^2}\big(|e^{ir\partial_x^2}g|^{p-1}e^{ir\partial_x^2}g\big)dr. \notag
\end{equation}

It remains to show that  $\calC_p$ is the sharp constant for the local inequality
\beq\label{ineq:G-N_type}
\int_0^1\|e^{ir\partial_x^2}f \| _{L^{p+1}}^{p+1}dr \leq \calC_p\|f\|_{L^2}^{\f{p+7}{2}} \|\partial_x f\|_{L^2}^{\f{p-5}{2}} .
\eeq
For the proof, we introduce the Weinstein functional $W_{p; I}$, with an interval $I\subset \R$, defined by
\[
 W_{p; I}(f)=\f{ \int_I\|e^{ir\partial_x^2}f \| _{L^{p+1}}^{p+1}dr }{\|f\|_{L^2}^{\f{p+7}{2}} \|\partial_x f\|_{L^2}^{\f{p-5}{2}}}\quad \textup{ for }f\in H^1(\R),
\]
and denote
\[
\calC_{p; I}= \sup_{f\in H^1(\R)} W_{p; I}(f).
\]
Then, it suffices to show $\mathcal{C}_{p;[0,1]}=\calC_p$. First, we note that
\begin{equation}\label{maximization scaling}
\mathcal{C}_{p;[0,1]}=\mathcal{C}_{p;[-\lambda^2,\lambda^2]}
\end{equation}
for any $\lambda>0$.
Indeed, if $f\in H^1(\R)$, then the function $e^{-\frac{i}{2}\partial_x^2}f_{\sqrt{2}\lambda,1}$, where $f_{\sqrt{2}\lambda,1}$ is defined in \eqref{eq:f_lambda,mu}, obeys
\[
\begin{aligned}
\int_0^1\|e^{ir\partial_x^2} e^{-\frac{i}{2}\partial_x^2}f_{\sqrt{2}\lambda,1}\| _{L^{p+1}}^{p+1}dr
&=\int_0^1\|(e^{i2\lambda^2(r-\frac{1}{2})\partial_x^2}f)(\sqrt{2}\lambda \cdot ) \| _{L^{p+1}}^{p+1}dr\\
&=\frac{1}{(\sqrt{2}\lambda)^{3}} \int_{-\lambda^2}^{\lambda^2} \|e^{ir\partial_x^2}f \| _{L^{p+1}}^{p+1}dr.
\end{aligned}
\]
Thus, together with \eqref{eq:scaling}, it follows that
\[\label{eq:equivalence1}
W_{p;[0,1]}(e^{-\frac{i}{2}\partial_x^2}f_{\sqrt{2}\lambda,1})=W_{p;[-\lambda^2,\lambda^2]}(f).
\]
Next, we claim that $\mathcal{C}_{p;[-\lambda^2,\lambda^2]}=\mathcal{C}_{p}$. Indeed, it is obvious that $\mathcal{C}_{p;[-\lambda^2,\lambda^2]}\leq\mathcal{C}_{p}$.
For the reverse inequality, given $\veps>0$, we choose a function $g\in H^1(\R)$ such that $\| g \|_{L^2}=\| \partial_x g\|_{L^2}=1$ and
\[
W_{p}(g)> \mathcal{C}_{p}-\veps.
\]
If $\lambda>0$ is large enough so that
$$\int_{-\lambda^2}^{\lambda^2} \|e^{ir\partial_x^2}g\| _{L^{p+1}}^{p+1}dr\geq \int_{\R} \|e^{ir\partial_x^2}g \| _{L^{p+1}}^{p+1}dr-\epsilon,$$
then
\[
\mathcal{C}_{p}<W_{p}(g)+\veps\leq W_{p; [-\lambda^2,\lambda^2]}(g)+2\veps \leq \mathcal{C}_{p;[-\lambda^2,\lambda^2]} +2\veps.
\]
Since $\veps>0$ is arbitrary, we conclude that $\mathcal{C}_p\leq\mathcal{C}_{p;[-\lambda^2,\lambda^2]}$.
\end{proof}

The following important norm quantities are independent of the possibly non-unique maximizer $Q$ which solves a certain Euler-Lagrange equation.
\begin{lemma}[Norm identities for $Q$]\label{norm identities}
Every maximizer $Q$ for the variational problem $\calC_p$, solving the Euler-Lagrange equation \eqref{eq:EL}, obeys
\beq\label{eq:Q}
\|Q\|_{L^2}^2=\left[\f{2(p+1)(p-5)}{(p+7)^2\calC_p}\right]^{\f{2}{p-1}}\cdot \left(\f{p+7}{p-5} \right)^{\f{1}{2}}
\quad \mbox{ and } \quad \|\partial_xQ\|^2_{L^2}=\f{p-5}{p+7}\|Q\|^2_{L^2} .
\eeq
Moreover,
\beq\label{eq:QandQ'}
\|e^{ir\partial_x^2}Q \|_{L_{r,x}^{p+1}}^{p+1}=\f{2(p+1)}{p+7}\|Q\|_{L^2}^2.
\eeq
\end{lemma}

\begin{proof}
Let $Q$ be a maximizer for $\calC_p$ solving \eqref{eq:EL} and define $g=\f{1}{\mu} Q(\f{1}{\lambda}\cdot)$ with
\begin{equation}\label{eq:scaling-lambda-mu}
\lambda=\frac{\|\partial_x Q\|_{L^2}}{\|Q\|_{L^2}}\quad \mbox{and}\quad
\mu=\|Q\|_{L^2}^{1/2}\|\partial_x Q\|_{L^2}^{1/2}.
\end{equation}
Since  $g$ is also a maximizer for $\calC_p$ with $\|g\|_{L^2}=\|\partial_x g\|_{L^2}=1$, it can be shown that $g$ solves
\[\label{eq:EL_g_2}
\f{p-5}{2} \partial_x^2 g-\f{p+7}{2} g
+\f{p+1}{\calC_p} \int_\R e^{-ir\partial_x^2}\big(|e^{ir\partial_x^2} g |^{p-1}e^{ir\partial_x^2} g \big)dr=0
\]
as in the proof of Proposition \ref{prop:maximizer}.
On the other hand,
 since $Q$ solves  \eqref{eq:EL}, $g$ solves
\[\label{eq:EL_mu}
\mu \lambda^2 \partial_x^2 g- \mu g + \mu^p\lambda^{-2}\int_\R e^{-ir\partial_x^2}\big(|e^{ir\partial_x^2} g |^{p-1}e^{ir\partial_x^2} g \big)dr=0.
\]
Combining the two equations, we have
\[
\left(\f{p-5}{2} \cdot \f{\calC_p}{p+1} - \f{\mu \lambda^2}{\mu^p\lambda^{-2}}\right)\partial_x^2 g
-\left(\f{p+7}{2} \cdot \f{\calC_p}{p+1} -  \f{\mu }{\mu^p\lambda^{-2}}\right)g=0.
\]
Since the operator $-\partial_x^2$ has no eigenvalue, the coefficients must be zero. Thus,
\[
\lambda=\left(\f{p-5}{p+7}\right)^{\f{1}{2}}\quad\mbox{and}\quad
\mu=\left[\f{2(p+1)(p-5)}{(p+7)^2\calC_p}\right]^{\f{1}{p-1}}
\]
which yields \eqref{eq:Q} by the choice of $\lambda$ and $\mu$ in \eqref{eq:scaling-lambda-mu}.
Moreover, the $L_x^2(\mathbb{R})$-inner product of the equation \eqref{eq:EL} by $Q$ yields
\[
\|\partial_x Q\|_{L^2}^2+\|Q\|_{L^2}^2-\|e^{ir\partial_x^2} Q\|_{L_{r,x}^{p+1}}^{p+1}=0
\]
which together with \eqref{eq:Q} deduces \eqref{eq:QandQ'}.
\end{proof}

\begin{remark}\label{rem:nonexistence}
For any $p\geq 9$, the variational problem $\mathcal{C}_{p;[0,1]}$ does not have a maximizer. 
Indeed, suppose that there exists a function $f \in H^1(\R)$ such that
\[
W_{p; [0,1]}(f) = \calC_{p;[0,1]}.
\]
Then,
 \[
 \begin{aligned}
\calC_{p;[-1,1]}  \geq W_ {p;[-1,1]}(f)
= \f{ \int_{-1}^1 \|e^{ir\partial_x^2}f\| _{L^{p+1}}^{p+1} dr}{\|f\|_{L^2}^{\f{p+7}{2}} \|\partial_x f\|_{L^2}^{\f{p-5}{2}}}
= \calC_{p;[0,1]} + W_{p;[-1,0]}(f)
\end{aligned}
 \]
which together with \eqref{maximization scaling} yields $\int_{-1}^0 \|e^{ir\partial_x^2}f\| _{L^{p+1}}^{p+1} dr=0$. Thus, we have $|e^{ir\partial_x^2}f|=0$ for almost all $r\in (-1,0)$ and $x\in \R$ and therefore
 \[
0=\int_{-1}^0 \|e^{ir\partial_x^2}f\|_{L^2}dr =\|f\|_{L^2}
 \]
which contradicts to $\|f\|_{L^2}>0$.
\end{remark}

\bigskip

\section{Proof of Theorem \ref{thm:globalvsblowup}}\label{sec:globalvsblowup}

In this section, we prove our main theorem for the global versus blowup dichotomy.

\begin{proof}[Proof of Theorem \ref{thm:globalvsblowup}]
Using the mass conservation and \eqref{ineq:G-N_type}, we have
\beq\label{eq:definition_f}
\begin{aligned}
   E[u(t)]M[u(t)]^\alpha &= \left(\f{1}{2}\|\partial_x u(t)\|^2_{L^2} - \f{1}{p+1} \int_0^1 \|e^{ir\partial_x^2} u(t)\|^{p+1}_{L^{p+1}}dr\right)\|u_0 \|_{L^2}^{2\alpha}\\
    &\geq  \f{1}{2}\left( \|\partial_x u(t)\|_{L^2}\|u_0 \|_{L^2}^\alpha \right)^2
    -  \f{\mathcal{C}_p}{p+1} \left(\|\partial_x u(t)\|_{L^2}\|u_0 \|_{L^2}^{\alpha} \right)^{\f{p-5}{2}}\\
    &=f(\|\partial_x u(t)\|_{L^2}\|u_0 \|^\alpha_{L^2})
\end{aligned}
\eeq
where $f:[0,\infty)\to\R$ is defined by
\[
f(x):=\f{1}{2}x^2- \f{\mathcal{C}_p}{p+1} x^{\f{p-5}{2}}.
\]
When $p>9$, $f$ has a unique local minimum at $x_0=0$ and the global maximum at
\[
    x_1 =\left[\f{(p-5)\mathcal{C}_p}{2(p+1)}\right]^{-\f{2}{p-9}}
    = \|\partial_xQ\|_{L^2}\|\|Q\|_{L^2}^\alpha
\]
where \eqref{eq: gradient times mass} is used in the second equality, see Figure \ref{figure} for the graph of $f$. Note also that $f(x_0)=0$ and
\[
f(x_1)=\left(\frac{1}{2}-\f{2}{p-5}\right)x_1^2=\frac{p-9}{2(p-5)}\|\partial_xQ\|^2_{L^2}\|Q\|_{L^2}^{2\alpha} =E_\infty[Q]M[Q]^\alpha
\]
where \eqref{eq:energyQ_invariant} is used in the last equality. Thus, it follows from \eqref{eq:definition_f}, the mass and the energy conservation laws, and the condition \eqref{ass:1} that
\beq\label{ineq:f}
f(\|\partial_x u(t)\|_{L^2}\|u_0 \|^\alpha_{L^2})\leq E[u(t)]M[u(t)]^\alpha=E[u_0 ]M[u_0 ]^{\alpha}<E_\infty[Q]M[Q]^\alpha=f(x_1).
\eeq%
%
%

Let us consider the first case such that condition \eqref{ass:global} holds, i.e., $\|\partial_xu_0\|_{L^2}\|u_0  \|_{L^2}^\alpha <x_1$. Then, by \eqref{ineq:f} and the continuity of $\|\partial_x u(t)\|_{L^2}$ in $t$,
we have $\|\partial_x u(t)\|_{L^2}\|u_0 \|^\alpha_{L^2}<x_1$ for all $t\in (-T_{\min},T_{\max})$, i.e., \eqref{conc:global} holds. Therefore, $\|u(t)\|_{H^1}$ stays bounded on $(-T_{\min},T_{\max})$ so that solution exists globally in time.

For the other case, we assume that condition  \eqref{ass:blowup} holds, i.e., $\|\partial_xu_0\|_{L^2}\|u_0  \|_{L^2}^\alpha > x_1$. Then, by the same argument as the above,  we have \eqref{conc:blowup} for all $t \in (-T_{\min},T_{\max})$. To show the last statement in Theorem \ref{thm:globalvsblowup}, suppose to the contrary that there exists an initial data $u_0  \in H^{1,1}(\mathbb{R})$ satisfying both \eqref{ass:1} and \eqref{ass:blowup} such that the solution $u(t)$ of \eqref{NLS} exists globally in time and therefore {\red $T_{\min}=\infty$}. Then, by density, we can choose a sequence $\{u_0  ^{(n)}\}_{n\in\N}\subset H^{3,3}(\mathbb{R})$ of initial data such that
\begin{equation}\label{density}
\|u_0  ^{(n)}-u_0  \|_{H^{1,1}}\to 0 \quad\mbox{as } n\to \infty.
\end{equation}
Note that
\[
E[u_0  ^{(n)}] \to E[u_0 ] \quad\mbox{ and }\quad M[u_0  ^{(n)}]\to M[u_0 ]
\]
as $n\to \infty$.
For each $n\in\N$, by Lemma \ref{Hkk LWP}, there exists the maximal solution $u^{(n)}\in C((-T_{\min}^{(n)}, T_{\max}^{(n)}); H^{3,3}(\R))$ to \eqref{NLS} with initial data $u_0  ^{(n)}$.
Observe that
\[
\begin{aligned}
\Phi_n(t)&= 16E[u_0  ^{(n)}]-\frac{4(p-9)}{p+1}\int_0^1\int_{\mathbb{R}}|e^{ir\partial_x^2}u^{(n)}(t,x)|^{p+1}dxdr\\
& = 4(p-5)E[u_0  ^{(n)}] -2(p-9)\|\partial_x u^{(n)}(t)\|_{L^2}^2,
\end{aligned}
\]
where $\Phi_n(t)$ is given in \eqref{eq:Phi} for $u^{(n)}$.
Multiplying both sides by $M[u_0 ]^\alpha$, we have
\[
\begin{aligned}
     \Phi_n(t)M[u_0 ]^\alpha  &= 4(p-5)E[u_0  ^{(n)}]M[u_0 ]^\alpha -2(p-9)\|\partial_x u^{(n)}(t)\|_{L^2}^2 \|u_0  \|_{L^2}^{2\alpha}.
    \end{aligned}
\]
It follows from \eqref{ass:1} that there exists $\delta>0$ such that $E[u_0 ]M[u_0 ]^\alpha <(1-\delta)E_\infty[Q]M[Q]^\alpha$.
Observe that $u_0  ^{(n)}$ satisfies the conditions \eqref{ass:1} and \eqref{ass:blowup} for sufficiently large $n$, since
$E[u_0  ^{(n)}]\to E[u_0 ]$  as $n\to \infty$. Thus,
 $\|\partial_x u^{(n)}(t)\|_{L^2}^2 \|u_0  \|_{L^2}^{2\alpha} > \|\partial_xQ\|_{L^2}^2 \|Q\|_{L^2}^{2\alpha}$
and, then, by \eqref{eq:energyQ_invariant}, we have
\beq\label{eq:negative}
\begin{aligned}
 \Phi_n(t)M[u_0 ]^\alpha
  &<  4(p-5)(1-\delta)E_\infty[Q]M[Q]^\alpha -2(p-9)\|\partial_xQ\|_{L^2}^2 \|Q\|_{L^2}^{2\alpha}\\
 & = -2\delta(p-9)\|\partial_xQ\|_{L^2}^2 \|Q\|_{L^2}^{2\alpha}<0
    \end{aligned}
\eeq
for sufficiently large $n$ and all $-T_{\min}^{(n)}<t<0$.
 If we denote  by $\mathcal{V}^{(n)}(t)$ the variance for $u^{(n)}$,
by Proposition \ref{prop:variance lower bound} and \eqref{eq:negative}, we have
\[
\begin{aligned}
\mathcal{V}^{(n)}(t) &\leq \mathcal{V}^{(n)}(0)+(\mathcal{V}_1^{(n)})'(0)t+\int_0^t\int_0^s \Phi_n(\tau )d\tau ds\\
& \leq  \mathcal{V}^{(n)}(0)+(\mathcal{V}_1^{(n)})'(0)t-\f{\delta (p-9)\|\partial_xQ\|_{L^2}^2 \|Q\|_{L^2}^{2\alpha}}{\|u_0  \|_{L^2}^{2\alpha}}\, t^2
\end{aligned}
\]
for sufficiently large $n$.
It immediately follows from \eqref{density} that
\[
\mathcal{V}^{(n)}(0) \to \int_\R x^2|u_0  (x)|^2 dx=:\mathcal{V}(0)\]
as $n\to \infty$.
Moreover,
$
(\mathcal{V}_1^{(n)})'(0) \to \mathcal{V}_1'(0)
$
as $n\to \infty$, where $\mathcal{V}_1'(0):=4\im \int_\R x \partial_xu_0(x) \overline{u_0  (x)}dx$, since
\[
\left|(\mathcal{V}_1^{(n)})'(0) - \mathcal{V}_1'(0) \right|
 \lesssim \|xu_0  ^{(n)}\|_{L^2}\|\partial_x(u_0  ^{(n)}-u_0  )\|_{L^2}+ \| \partial_xu_0 \|_{L^2}\|x(u_0  ^{(n)}-u_0  )\|_{L^2}.
\]
Therefore, we obtain
\beq\label{ineq:V_n}
\mathcal{V}^{(n)}(t) \lesssim \mathcal{V}(0)|\mathcal{V}_1'(0)|t-\f{\delta (p-9)\|Q'\|_{L^2}^2 \|Q\|_{L^2}^{2\alpha}}{\|u_0  \|_{L^2}^{2\alpha}}\, t^2
\eeq
for sufficiently large $n$.
This yields the existence of a positive $T$ such that the right-hand side of \eqref{ineq:V_n} is negative for all $t<-T$ since it is a second-degree polynomial and the coefficient of $t^2$ is negative. 
However, since $\mathcal{V}^{(n)}(t) \geq 0$ for all $-T_{\min}^{(n)}<t<0$, we observe that  for sufficiently large $n$,  $u^{(n)}(t)$ cannot exist beyond the time $T$ and therefore $ -T_{\min}^{(n)}\geq -T$. 
We also note that by the persistence of regularity (Lemma \ref{lemma: persistence of regularity}), $\|u^{(n)}(t)\|_{H^1}$ blows up as $t\downarrow -T_{\textup{min}}^{(n)}\geq -T$ for sufficiently large $n$.

On the other hand, by the hypothesis, $u$ exists in $C([-T,0]; H^1(\mathbb{R}))$ and, moreover, Lemma \ref{lemma: continuity of data-to-solution map} implies that $u^{(n)}\in C([-T,0]; H^1(\mathbb{R}))$ for sufficiently large $n$. Therefore, we deduce a contradiction.
\end{proof}

\bigskip

\noindent

\textbf{Acknowledgements: }
The authors are supported by the National Research Foundation of Korea (NRF) grants funded by the Korean government (MSIT) NRF-2020R1A2C1A01010735, (MSIT) RS-2023-00208824, (MSIT) RS-2023-00219980 and (MOE) NRF-2021R1I1A1A01045900.

\end{document}